\documentclass[11pt]{amsart}

\usepackage[T1]{fontenc}
\usepackage[utf8]{inputenc}
\usepackage[english]{babel}
\usepackage{lmodern,microtype}
\usepackage{amsmath,amssymb,amsthm,mathtools}
\usepackage{mathrsfs}
\usepackage[margin=30mm]{geometry}
\usepackage[colorlinks=true,linkcolor=blue,citecolor=blue,urlcolor=blue]{hyperref}
\usepackage[nameinlink,noabbrev]{cleveref}

\newtheorem{theorem}{Theorem}[section]
\newtheorem{proposition}[theorem]{Proposition}
\newtheorem{lemma}[theorem]{Lemma}

\newtheorem{conjecture}[theorem]{Conjecture}
\theoremstyle{definition}

\newtheorem{example}[theorem]{Example}
\theoremstyle{remark}

\def\<{\langle}
\def\>{\rangle}

\newcommand{\QQ}{\mathbb Q}
\newcommand{\NN}{\mathbb N}
\newcommand{\Sym}{\Lambda}
\newcommand{\ellp}{\ell}
\newcommand{\whp}{\widehat p}
\newcommand{\Ch}{\operatorname{Ch}}

\newcommand{\cB}{\mathsf B}
\newcommand{\cT}{\mathcal T}
\newcommand{\Sn}{\mathfrak S}

\title{An operator $b$-positivity conjecture \\ for Jack characters}
\author{Jean-Yves Thibon}
\address{Laboratoire d’Informatique Gaspard-Monge,
Université Gustave Eiffel,
CNRS, ESIEE Paris,
F-77454 Marne-la-Vallée, France}
\email[Jean-Yves Thibon]{jean-yves.thibon@univ-eiffel.fr} 


\date{}

\begin{document}
\maketitle

\begin{abstract}
Goulden and Jackson defined a Jack deformation $a_{\mu\nu}^\lambda(\alpha)$ of the
connection coefficients of the symmetric group and conjectured that these were polynomials
in $b=\alpha-1$ with nonnegative integer coefficients. This conjecture was extended to
the structure constants of (stable) Jack characters by Śniady.  
We formulate a stronger operator conjecture: the differential operators $\Delta_\mu(\alpha)$
that multiply by stable Jack characters should have coefficients in $\NN[b]$ in the normally
ordered Heisenberg basis $p_ID_J$. We prove integrality of the normally
ordered coefficients and give further positivity results for single-part
indices. We also provide numerical evidence, as well as theoretical consequences
of the conjecture which can be rigorously proved. Namely, we exhibit two families of $\alpha$-shifted
symmetric functions whose associated differential operators are $\NN[b]$ positive 
simultaneously on the bases $p_ID_J$ and $\Delta_\mu(\alpha)$. These operators correspond to the Boolean
	cumulants of the Jack transition measure and to the Cartan generators of the affine Yangian  $Y(\widehat{\mathfrak gl}_1)$.
\end{abstract}


\section{From class products to stable operators}

\subsection{The classical model}
For a partition $\mu\vdash n$, let $K_\mu$ denote the sum of the
permutations of cycle type $\mu$ in the center of $\QQ[S_n]$.  The
connection coefficients are defined by
\[
 K_\mu K_\nu=\sum_{\lambda\vdash n}
 a_{\mu\nu}^{\lambda}K_\lambda.
\]
The center $Z_n$ of the symmetric group algebra $\QQ[\Sn_n]$ can be identified with the space
$\Lambda_n$ of symmetric functions of degree $n$ via the linear map $Z(\sigma)=p_\mu$ sending a permutation
$\sigma$ of cycle type $\mu$ to the power sum product $p_\mu$. 
The product of $Z_n$ induces a product $\times$ on $\Lambda_n$ called the class product. On Schur functions, it
is given by
\[ 
s_\lambda\times s_\mu = \delta_{\lambda\mu}\frac1{f^\lambda}s_\lambda
\]
where $f^\lambda$ is the dimension of the irreducible representation indexed by $\lambda$.
Thus, the class product is dual to the coproduct
\[
	\Gamma(s_\lambda)=\frac1{f^\lambda}s_\lambda\otimes s_\lambda
\]
and
\[
\Gamma(p_\lambda)=\frac1{n!}\sum_{\mu,\nu}a_{\mu\nu}^\lambda p_\mu(X)p_\nu(Y)
=\<p_\lambda, \Phi\>_Z,
\]
where
\[
	\Phi(X,Y,Z)=\sum_\lambda \frac1{f^\lambda}s_\lambda(X)s_\lambda(Y)s_\lambda(Z).
\]
Ivanov and Kerov \cite{IvanovKerov} have proved that the structure constants of normalized conjugacy
classes ${\mathsf a}_{\mu;n}$ such that $Z({\mathsf a}_{\mu;n})=(n)_r p_\mu p_1^{n-r}$ are independent of $n$.
In terms of the class product, extended to the degree completion of $\Lambda$, this means that the formal series
\begin{equation}\label{eq:stable-series}
 \whp_\mu=\frac{|\mu|!\,p_\mu}{(1-p_1)^{|\mu|+1}}.
\end{equation}
are stable under $\times$. Indeed, there exist nonnegative integers $d_{\mu\nu}^\lambda$ such that
\begin{equation}\label{eq:stable-d}
 \whp_\mu\times\whp_\nu
 =\sum_\lambda d_{\mu\nu}^{\lambda}\whp_\lambda.
\end{equation}
This can be deduced from the existence of differential operators implementing the class product
by $\whp_\mu$
\begin{equation}\label{eq:Delta-classical}
 \Delta_\mu(f)=\whp_\mu\times f
\end{equation}
Indeed, these operators are nonnegative sums of terms $p_ID_J$ whose action on a $\whp_\nu$ produces only 
terms of the same type. Their existence has been conjectured by Katriel \cite{Katriel} and first proved
by Goupil-Poulalhon-Schaeffer \cite{GoupilPoulalhonSchaeffer}. A different proof in terms of matrix integrals
due to Lascoux and the author \cite{LascouxThibon}, has the advantage of remaining valid for the values
$\alpha=1,2,\frac12$ of the deformed class product, and sheds some light on the 
origin of matchings in the matching Jack conjecture recalled below.

The operators $\Delta_\mu$ are diagonal on the Schur basis,and their eigenvalues are the shifted
power-sums, also known as (stable) central characters
\begin{equation}
	\Delta_\mu s_\lambda = p_\mu^\#(\lambda)s_\lambda = \Ch_\mu(\lambda)s_\lambda.
\end{equation}
The knowledge of the differential operators is sufficient to reconstruct from scratch the stable series,
hence the Ivanov-Kerov algebra, and the theory of shifted symmetric functions \cite{ThibonStableJack}.

\subsection{Jack deformation and shifted Jack functions}

Let $J_\lambda^{(\alpha)}$ be the integral Jack functions, orthogonal for
\[
 \langle p_\mu,p_\nu\rangle_\alpha
 =\delta_{\mu\nu}z_\mu\alpha^{\ellp(\mu)},
 \qquad b=\alpha-1.
\]
Goulden and Jackson define $a_{\mu\nu}^{\lambda}(\alpha)$ by
\begin{equation}\label{eq:GJ-series}
 \sum_\lambda\frac{J_\lambda^{(\alpha)}(X)
                         J_\lambda^{(\alpha)}(Y)
                         J_\lambda^{(\alpha)}(Z)}
                        {j_\lambda^{(\alpha)}}t^{|\lambda|}
 =1+\sum_{n\geq1}t^n
   \sum_{\mu,\nu,\lambda\vdash n}
   \frac{a_{\mu\nu}^{\lambda}(\alpha)}
        {z_\lambda\alpha^{\ellp(\lambda)}}
        p_\mu(X)p_\nu(Y)p_\lambda(Z).
\end{equation}
This defines a deformation $\times_\alpha$ of the class product.  
If we choose to interpret the coefficients $a_{\mu\nu}^\lambda(\alpha)$ as deformed structure constants
of the conjugacy classes, this imposes to define $\times_\alpha$ as dual to
the diagonal Jack coproduct
\[
 \Gamma_\alpha(J_\lambda^{(\alpha)})
 =\frac1{n!}J_\lambda^{(\alpha)}\otimes J_\lambda^{(\alpha)},
 \qquad |\lambda|=n,
\]
for the usual Hall scalar product:
\begin{equation}\label{eq:Jack-product}
 \langle f\times_\alpha g,h\rangle_1
 =\langle f\otimes g,\Gamma_\alpha(h)\rangle_{1\otimes1}.
\end{equation}
It might seem more logical to dualize the coproduct with respect to the Jack scalar
product. This would only introduce normalizations by powers of $\alpha$ in the formulas.
Our choice has the advantage that the series $\whp_\mu$ remain stable under the deformed product
without any modification. It has the drawback that its forces to dualize all operators 
taken from the integrable systems literature. For example, the $\Delta_2(\alpha)$ is the Hall adjoint
of the Jack Laplace-Beltrami operator.

We denote by $\Delta_\mu^\perp$  the adjoint of $\Delta_\mu$ with respect to the Hall scalar product.

The shifted Jack functions of Okounkov--Olshanski and their later
developments provide eigenvalue coordinates for this product
\cite{OkounkovOlshanski,BenDaliDolega}.  For $\mu\vdash m$, use the same
stable series $\whp_\mu$ as in \eqref{eq:stable-series}, and define
\begin{equation}\label{eq:Jack-Delta}
 \Delta_\mu(\alpha)(f)=\whp_\mu\times_\alpha f.
\end{equation}
Their eigenvalues on the modified Jack functions $J'_\lambda(X)=J^{(\alpha)}_\lambda[\alpha X]$,
where $p_n[\alpha X]=\alpha p_n(X)$, are the normalized Jack characters
$\Ch_\mu^{(\alpha)}(\lambda)$.  
The stable structure constants
$d_{\mu\nu}^{\lambda}(\alpha)$ are specified by
\begin{equation}\label{eq:Jack-d}
 \whp_\mu\times_\alpha\whp_\nu
 =\sum_\lambda d_{\mu\nu}^{\lambda}(\alpha)\whp_\lambda.
\end{equation}
Equivalently, they are the structure constants for multiplication of the
corresponding stable Jack characters.

Thus, the $\times_\alpha$ stability of the series $\whp_\mu$ follows
from the stability under product of Jack shifted symmetric functions.
The existence of the correponding differential operators is implied by that of
the Nazarov-Sklyanin operators. The $\Delta_\mu(\alpha)$ can be computed recursively
from the operators $A^{(k)}$ of \cite{NazarovSklyanin} by the Newton formulas
(see \cite{ThibonStableJack} for details and examples).

\section{Positivity conjectures}

\begin{conjecture}[Matching--Jack, \cite{GouldenJackson1996}]\label{conj:matching-jack}
For partitions $\lambda,\mu,\nu$ of the same size,
\[
 a_{\mu\nu}^{\lambda}(1+b)\in\NN[b].
\]
\end{conjecture}
The original matching formulation asks for more than this positivity: for
a family $\mathcal M_{\mu\nu}^{\lambda}$ of matchings associated with
$(\lambda,\mu,\nu)$, it asks for a statistic $\operatorname{stat}$ such
that
\[
 a_{\mu\nu}^{\lambda}(1+b)
   =\sum_{M\in\mathcal M_{\mu\nu}^{\lambda}}
        b^{\operatorname{stat}(M)}.
\]
The statistic is expected to measure non-bipartiteness: it vanishes exactly
on the bipartite matchings.  The operator conjecture below would imply the
coefficientwise positivity, but by itself supplies no such matching model
or statistic.

At $\alpha=1$ and $\alpha=2$, the complex and real Gaussian matrix
integrals give Wick expansions in which the relevant pairings and the
bipartite/non-bipartite distinction are visible.  In particular, the
$\alpha=2$ expansion provides the required statistic  on the
real-matrix side.  The quaternionic specialization $\alpha=1/2$ also has a
Wick expansion, but its contraction weights include signs and factors that
do not give the same direct nonnegative statistic.  These integral formulas
explain the specializations but does not suggest the matching formula for
general $b$ \cite{GouldenJackson1996,ThibonStableJack}.

\begin{conjecture}[Śniady, \cite{Sniedy}]\label{conj:sniady}
For arbitrary partitions $\lambda,\mu,\nu$,
\[
 d_{\mu\nu}^{\lambda}(1+b)\in\NN[b].
\]
The equal-size coefficients in \eqref{eq:GJ-series} are the corresponding
homogeneous specialization, with the normalization fixed there.
\end{conjecture}

Polynomiality is known in this setting; coefficientwise positivity is the
open part of these conjectures \cite{DolegaFerayBCumulants,Sniedy}.

Let $D_r=r\,\partial/\partial p_r$, so that
\begin{equation}\label{eq:Heisenberg}
 [D_r,p_s]=r\delta_{rs},\qquad [p_r,p_s]=[D_r,D_s]=0.
\end{equation}
For partitions $I=(i_1,\ldots)$ and $J=(j_1,\ldots)$, put
$p_I=\prod_a p_{i_a}$ and $D_J=\prod_bD_{j_b}$.  Normal order means that
all multiplications $p_r$ are placed to the left of all derivatives $D_s$.

\begin{conjecture}[Operator $b$-positivity]\label{conj:operator}
For every partition $\mu$, write the stable multiplier in normal order as
\begin{equation}\label{eq:normal-order-Delta}
 \Delta_\mu(1+b)=\sum_{I,J}c_{I,J}^{\mu}(b)p_ID_J.
\end{equation}
Then
\[
 c_{I,J}^{\mu}(b)\in\NN[b]\qquad\text{for all }I,J,\mu.
\]
\end{conjecture}

This is stronger than the two structure-constant conjectures.  The following
calculation explains the implication. 

\begin{lemma}[Positive action of balanced Heisenberg monomials]
\label{lem:balanced-action}
Suppose $|I|=|J|$, let $i_r=m_r(I)$, $j_r=m_r(J)$, and let
$n_r=m_r(\nu)$.  If $j_r\leq n_r$ for all $r\geq2$, then
\begin{align}\label{eq:balanced-action}
 p_ID_J\whp_\nu
 ={}&\left(\prod_{r\geq2}r^{j_r}(n_r)_{\underline{j_r}}\right)
 \sum_{s=0}^{\min(j_1,n_1)}
 \binom{j_1}{s}(n_1)_{\underline{s}}\notag\\
 &\hspace{18mm}\times
 \whp_{\, (\nu\setminus(J_{\geq2}\cup1^s))\cup I}.
\end{align}
If any $j_r>n_r$ for $r\geq2$, the result is zero.  In particular,
all coefficients in this action are nonnegative integers.
\end{lemma}

\begin{proof}
Use $\whp_\nu=|\nu|!p_\nu(1-p_1)^{-|\nu|-1}$.  For $r\geq2$,
$D_r^{j_r}$ acts on the numerator and contributes
$r^{j_r}(n_r)_{\underline{j_r}}$.  In the $p_1$ variable, choose $s$ of
the $j_1$ derivatives to act on the numerator; the remaining derivatives
act on $(1-p_1)^{-|\nu|-1}$.  Leibniz's rule gives
\[
 \binom{j_1}{s}(n_1)_{\underline{s}}
 (|\nu|+1)^{\overline{j_1-s}}.
\]
The resulting index has size $|\nu|+j_1-s$.  The ratio of factorials in
the stable-series normalization cancels the rising factorial, leaving
\eqref{eq:balanced-action}.
\end{proof}

Each $\Delta_\mu$ preserves weighted degree, so every monomial $p_ID_J$
in \eqref{eq:normal-order-Delta} has $|I|=|J|$.  The lemma therefore shows
that Conjecture~\ref{conj:operator} implies coefficientwise positivity of
the stable products \eqref{eq:Jack-d}.  In the dehomogenized three-Jack
kernel, these are the extended-degree structure constants of stable Jack
characters; when $|\lambda|=|\mu|=|\nu|$, they agree with the
Goulden--Jackson coefficients $a_{\mu\nu}^{\lambda}(\alpha)$ in the
normalization of \eqref{eq:GJ-series}
\cite[Proposition~6.2]{BenDaliSuperNabla}.  Thus the operator conjecture
implies Śniady positivity and the coefficientwise positivity part of
Matching--Jack. 

\subsection*{Two positive families studied below}

The evidence developed in the rest of the note comes from two families of
operators.  The Boolean cumulants $\cB_n^{(\alpha)}$ are the coefficients
of the inverse Nazarov--Sklyanin transfer operator.  Their first-return Lax
paths give normally ordered expansions with coefficients in $\NN[b]$;
their ribbon-path formula gives a second $\NN[b]$-positive expansion in
the stable-character operators $\Delta_\mu^\perp$.  
The second family, given by
$H_q=[D_1,\operatorname{ad}_{\mathcal D_\alpha}^{q}(p_1)]$, where
$\mathcal D_\alpha$ is the first Jack Laplace--Beltrami Hamiltonian, are
proportional to the Cartan generators $\psi_q$ of the affine Yangian, and are
constructed directly in the Weyl algebra.  Wick ordering their path formula
proves positivity in $p_ID_J$, and a separate insertion argument proves
positivity in the $\Delta_\mu^\perp$ basis.

The normal forms already show the pattern.  With
$E=\sum_{r\geq1}p_rD_r$, the first Boolean modes are
\begin{align*}
 \cB_2^{(\alpha)}={}&\alpha\sum_{r\geq1}p_rD_r = \alpha E,\\
 \cB_3^{(\alpha)}={}&\alpha\sum_{a,b\geq1}p_ap_bD_{a+b}
       +\alpha^2\sum_{a,b\geq1}p_{a+b}D_aD_b
       +\alpha b\sum_{r\geq1}r p_rD_r.
\end{align*}
Each summand is already a normally ordered monomial sum, and after
$\alpha=1+b$ its scalar coefficient is in $\NN[b]$.  On the Cartan side,
$A_1=\alpha\sum_{h\geq1}p_{h+1}D_h$ and $H_2=2\alpha E$. 
For higher
modes, the path words are ordered first and Wick contractions use only
$[\alpha D_r,p_s]=r(1+b)\delta_{rs}$.  The detailed path formulas and the
corresponding stable-character expansions follow in Sections~4 and~5.

\section{Exact computations and numerical evidence}

We have expanded the operators $\Delta_\mu(1+b)$ in the normal-ordered
basis and checked the resulting coefficients.  The present
verification covers $|\mu|\leq5$ and $|I|=|J|\leq14$, for a total of
$31\,910$ coefficients.  The arithmetic is in $\QQ[b]$; every coefficient
tested belongs to $\NN[b]$.  The explicit formulas for $\Delta_3$ and
$\Delta_4$ and the stable-series construction are recorded in
\cite{ThibonStableJack}.

The first cumulants are
\begin{align}
 \cB_2^{(\alpha)}={}&\alpha E,\notag\\
 \cB_3^{(\alpha)}={}&\alpha(\Delta_2^\perp+bE),\notag\\
 \cB_4^{(\alpha)}={}&\alpha\bigl(\Delta_3^\perp
 +3b\Delta_2^\perp+\alpha(E)_2+b^2E\bigr),\label{eq:first-Boolean}\\
 \cB_5^{(\alpha)}={}&\alpha\bigl(\Delta_4^\perp
 +6b\Delta_3^\perp+3\alpha\Delta_{21}^\perp
 +(7b^2+b+1)\Delta_2^\perp
 +4b\alpha(E)_2+b^3E\bigr),\notag
\end{align}
and the first Cartan modes are
character basis:
\begin{align*}
 H_3={}&(3+3b)\Delta_2^\perp+(2b+2b^2)\Delta_1^\perp,\\
 H_4={}&(4+4b)\Delta_3^\perp+6(1+2b+b^2)\Delta_{11}^\perp
       +(9b+9b^2)\Delta_2^\perp
       +2(1+2b+2b^2+b^3)\Delta_1^\perp,\\
 H_5={}&(5+5b)\Delta_4^\perp+(24b+24b^2)\Delta_3^\perp
       +(20+40b+20b^2)\Delta_{21}^\perp
       +(18b+36b^2+18b^3)\Delta_{11}^\perp\\
 &\quad +(15+30b+36b^2+21b^3)\Delta_2^\perp
       +(4b+8b^2+6b^3+2b^4)\Delta_1^\perp.
\end{align*}

\subsection{Additional positivity results for a single part}
\label{subsec:single-part-checks}

Write $T_n=\Delta_{(n)}^\perp$ and, for $|I|=|J|=N\geq n\geq2$, put
\[
 F_{n,N}(I,J)=\alpha^{1-\ellp(J)}[p_ID_J]T_n,
 \qquad d=n+1-\ellp(I)-\ellp(J).
\]
Here $N$ is the weight of the normal monomial, while $n$ is the
character index. The integrality of the normalized coefficient
$F_{n,N}(I,J)$ follows from Section~\ref{subsec:integrality}. In the
notation of that section, the connected three-Jack series gives
\begin{equation}\label{eq:single-part-connected}
 F_{n,N}(I,J)=n[t^N p_n[A]p_1[A]^{N-n}p_I[X]p_J[Y]]
 \bigl(\alpha\log\tau_\alpha(t;A,X,Y)\bigr).
\end{equation}
Indeed, in the logarithm of the quotient in
\eqref{eq:normalized-symbol-kernel}, the term linear in $p_n[Z]$ is
the corresponding linear term of the quotient itself. The substitution
$p_n[A]=\alpha p_n[Z]$ for $n\geq2$ gives the factor $\alpha$, and
$z_{(n)}=n$ gives the remaining factor. The degree bound of
Do\l{}\k{e}ga--F\'eray gives $\deg_b F_{n,N}(I,J)\leq d$
\cite[Theorem~1.2]{DolegaFerayBCumulants}. In particular the coefficient
vanishes for $d<0$.

The following cases hold for all the indicated weights.
\begin{itemize}
\item For every $n\geq2$, the coefficients at $N=n$ are in $\NN[b]$.
Indeed, they are obtained by taking the term linear in $u$ in the
connected series with $p_r[A]=u$ for every $r$. At weight $n$, this
selects the single partition $(n)$. The connected marginal positivity
theorem of Chapuy--Do\l{}\k{e}ga applies
\cite[Theorem~5.10]{ChapuyDolega}.
\item For every $n\geq2$ and every $N\geq n$, the coefficients with
$\ellp(I)=1$ or $\ellp(J)=1$ are in $\NN[b]$. For example, if $I=(N)$,
exchanging $A$ and $X$ in \eqref{eq:single-part-connected} gives
\[
 F_{n,N}((N),J)=\frac{n}{N}F_{N,N}((n,1^{N-n}),J).
\]
The ratio is positive, and integrality of the left-hand side gives the
claim.
\item For every $n\geq2$ and every $N\geq n$, the coefficients of defect
$d\leq4$ are in $\NN[b]$. In these sectors, the duality
$F(\alpha)=(-1)^d\alpha^d F(1/\alpha)$ and the specializations
$\alpha=0,1,2$ determine $F=F_{n,N}(I,J)$. The verification uses
insertion series on the connected bipartite incidence graphs with
$n$ edges. At each graph, the three specialized weights count
orientable rotations, orientable one-face systems, and signed one-face
systems, respectively, with the sign $(-1)^d$ at $\alpha=0$.
After contraction of a spanning tree, the signed one-face condition is
invertibility of a $d\times d$ alternating intersection matrix over
$\mathbb F_2$ with a variable diagonal. Enumerating the diagonals for
$d\leq4$ gives nonnegative weights in the basis
$\alpha^j b^{d-2j}$. The insertion series have nonnegative coefficients.
\end{itemize}

Further exact verifications give $b$-positivity of
$\Delta_{(n)}$ at all weights for $1\leq n\leq9$. For $n=10$
they cover all sectors except $\ellp(I)=\ellp(J)=3$. In particular,
the $2\times2$ sector has an all-weight certificate. The homogeneous
insertion representative of degree $10$ was checked at weights
$10,\ldots,15$ and at the eight values $\alpha=1,\ldots,8$.
The weight bound follows from the following observation. For
$t_{ij}=vU_iV_j$ and $g_{ij}=t_{ij}/(1-t_{ij})$, the leading
specialization of a homogeneous polynomial in four variables is
obtained by replacing $g_{ij}$ by $t_{ij}$. The kernel of this
leading specialization is generated by the $2\times2$ minor $Q$.
Since $Q(g)$ has $v$-valuation $3$, a homogeneous polynomial
of degree $10$ divisible by $Q^k$ and by no higher power has
$v$-valuation $10+k\leq15$. The degree bound in $b$ in this
sector is $7$, which explains the eight specializations.

For the remaining $3\times3$ sector of $\Delta_{(10)}$, exact
interpolation in the three-dimensional duality space at defect five
confirms positivity for $10\leq N\leq14$. This check comprises
$410$ profile pairs up to interchange of $I$ and $J$ and establishes
positivity only in that weight range. In defect five, duality gives
$F=Rb^5+M\alpha b^3+Q\alpha^2b$, where $R=-F(0)$ and $Q=F'(1)$
when $F$ is viewed as a polynomial in $\alpha$. Since $R\geq0$, the
two inequalities $M\geq0$ and $Q\geq0$ are equivalent to
$b$-positivity: $M$ and $Q$ are the coefficients of $b^4$ and $b$,
respectively, and all other coefficients are nonnegative sums of
$R,M,Q$. All statements about $F_{n,N}(I,J)$ transfer to $T_n$ by multiplication by
$\alpha^{\ellp(J)-1}$ and then to $\Delta_{(n)}$ by Hall adjunction.

\section{Boolean cumulants and the Nazarov--Sklyanin operators}

The first positive family comes from the Nazarov--Sklyanin transfer
operator.  
Put $b=\alpha-1$ and define the Lax matrix, indexed by nonnegative integers,
by
\begin{equation}\label{eq:Lax}
 L_{ij}=\begin{cases}
 p_{j-i},&j>i,\\
 ib,&j=i,\\
 \alpha D_{i-j},&j<i.
 \end{cases}
\end{equation}
On $\Sym[w]=\bigoplus_{i\geq0}\Sym w^i$, let $\iota_0(f)=f$ in the
vacuum component and let $\pi_0$ project onto that component.  Set
\begin{equation}\label{eq:transfer}
 \cT(u)=\pi_0(u-L)^{-1}\iota_0.
\end{equation}
Define the diagonal operators $\cB_n^{(\alpha)}$ by
\begin{equation}\label{eq:Boolean}
 \cT(u)^{-1}=u\left(1-\sum_{n\geq1}\cB_n^{(\alpha)}u^{-n}\right).
\end{equation}

\begin{proposition}[First-return formula]\label{prop:first-return}
For $n\geq2$,
\[
 \cB_n^{(\alpha)}=L_{0+}L_{++}^{n-2}L_{+0},
\]
where $L_{++}$ is the restriction to the positive auxiliary states.
Thus $\cB_n^{(\alpha)}$ is the sum of weights of Lax paths that leave
height zero at the first step, remain strictly positive, and return to zero
at the last step.
\end{proposition}

\begin{proof}
Write $L$ in vacuum/positive blocks:
\[
 L=\begin{pmatrix}0&L_{0+}\\L_{+0}&L_{++}\end{pmatrix}.
\]
Solving $(u-L)(x_0,x_+)=(f,0)$ gives
$x_+=(u-L_{++})^{-1}L_{+0}x_0$ and hence
\[
 \cT(u)^{-1}=u-L_{0+}(u-L_{++})^{-1}L_{+0}.
\]
Expanding the resolvent in powers of $u^{-1}$ and comparing with
\eqref{eq:Boolean} gives the formula.  The positive intermediate block
forces the path to avoid zero until its final return.
\end{proof}

The ribbon-path formula of
Ben Dali--Do\l{}\k{e}ga \cite{BenDaliDolega} identifies the general
stable-character expansion 
\begin{equation}\label{eq:B-ribbon}
 \cB_n^{(\alpha)}
 =\sum_{R\in\mathcal L_n}
 \omega(R)\alpha^{\ellp(\mu(R))+a(R)}b^{h(R)}
 \Delta_{\mu(R)}^\perp.
\end{equation}
Here $R$ is a first-return Łukasiewicz excursion decorated by disjoint
pairs.  A pair of degree $r$ joins a downward step of size $r$ to a later
upward step of the same size; pairs may be nested or cross.  Every unpaired
upward step has size one, equivalently every upward step of size at least
two is paired.  The partition $\mu(R)$ records the unpaired downward steps,
$a(R)$ counts pairs, $h_i(R)$ counts horizontal steps at height $i$, and
$h(R)=\sum_i h_i(R)$.  The weight is
$\omega(R)=\prod_{P\text{ pair}}\deg(P)\prod_{i\geq1}i^{h_i(R)}$.
All these weights are nonnegative integers, so
\eqref{eq:B-ribbon} is coefficientwise positive in $b$.

\begin{proposition}[Boolean multiplication on stable series]
\label{prop:B-multiplication}
For partitions $\mu,\nu$ with every part of $\mu$ at least $2$, there are coefficients
$c_{\mu,\nu}^{\lambda}(b)\in\NN[b]$ such that
\[
 \mathsf B_\mu\times_\alpha\whp_\nu
 =\sum_\lambda c_{\mu,\nu}^{\lambda}(b)\whp_\lambda,
\]
where, for each $n\geq2$, $\mathsf B_{(n)}$ is the stable series whose multiplier is
$(\cB_n^{(\alpha)})^\perp$, and
$\mathsf B_\mu=\mathsf B_{(\mu_1)}\times_\alpha\cdots\times_\alpha
\mathsf B_{(\mu_{\ellp(\mu)})}$.
\end{proposition}

\begin{proof}
The first-return formula expands $\cB_n^{(\alpha)}$ into Lax path words
whose factors are $p_r$, $\alpha D_r$, and nonnegative multiples of $b$.
Normal ordering uses only $[D_r,p_s]=r\delta_{rs}$, so every coefficient
is in $\NN[b]$.  Hall adjunction sends $p_ID_J$ to $p_JD_I$ and preserves
this positivity.  Apply Lemma~\ref{lem:balanced-action} term by term to
$(\cB_\mu^{(\alpha)})^\perp\whp_\nu$.  Each term becomes a nonnegative
integer combination of stable series, proving the assertion.
\end{proof}

\section{Cartan modes constructed from the first Jack Hamiltonian}

The second family is constructed directly in the Weyl algebra.  It does not
use the transfer operator to establish existence; the transfer operator
will enter only in the proof of stable-character positivity.

\subsection{Construction, diagonality, and shifted coordinates}

Let $\mathcal D_\alpha$ be half the Jack Laplace--Beltrami operator:
\begin{equation}\label{eq:LB}
 2\mathcal D_\alpha=
 \alpha\sum_{i,j\geq1}p_{i+j}D_iD_j
 +\sum_{i,j\geq1}p_ip_jD_{i+j}
 +b\sum_{r\geq1}(r-1)p_rD_r.
\end{equation}
For a box $\square=(i,j)$ set
$c_\alpha(\square)=\alpha(j-1)-(i-1)$.  Define
\begin{equation}\label{eq:A-H}
 A_q=\operatorname{ad}_{\mathcal D_\alpha}^{q}(p_1),
 \qquad H_q=[D_1,A_q].
\end{equation}
The operator $A_q$ adds a box with its one-box Jack--Pieri coefficient
multiplied by $c_\alpha(\square)^q$.

\begin{lemma}[Direct diagonality]\label{lem:H-diag}
Each $H_q$ is diagonal in the Jack basis.  If $F_q(\lambda)$ is its
eigenvalue, then
\begin{align}\label{eq:Fq-Pieri}
 F_q(\lambda)={}&\sum_{\nu\gtrdot\lambda}
 c_\alpha(\nu/\lambda)^q u_{\lambda\nu}d_{\nu\lambda}\notag\\
 &-\sum_{\rho\lessdot\lambda}
 c_\alpha(\lambda/\rho)^q d_{\lambda\rho}u_{\rho\lambda},
\end{align}
where $u$ and $d$ are the one-box coefficients for $p_1$ and $D_1$.
\end{lemma}

\begin{proof}
For distinct $\lambda,\mu$ of the same size, a coefficient of
$D_1A_q-A_qD_1$ can occur only through a common cover and a common lower
cover.  If they exist, they are uniquely $\lambda\cup\mu$ and
$\lambda\cap\mu$.  The added box is the same along both routes, so both
terms carry the same content power.  The remaining difference is the
off-diagonal coefficient of $[D_1,p_1]=1$, hence vanishes.  The diagonal
coefficient is \eqref{eq:Fq-Pieri}.
\end{proof}

From the explicit formulas for the coefficients $u$ and $d$, one can prove the 
generating function
\begin{equation}\label{eq:Cartan-charge}
 1+b\sum_{q\geq0}F_q(\lambda)z^{-q-1}
 =\left(1+\frac bz\right)
   \prod_{\square\in\lambda}
   \prod_{s\in\{\alpha,-1,-b\}}
   \frac{z-c_\alpha(\square)+s}{z-c_\alpha(\square)-s}.
\end{equation}
It allows to identify $H_q$  with the Cartan modes $\alpha\psi_q$ of the affine
Yangian of $\mathfrak{gl}_1$; the three parameters are
$(s_1,s_2,s_3)=(\alpha,-1,-b)$ \cite{LiAffineYangian}.  It also agrees
with the Nazarov--Sklyanin transfer quotient
\begin{equation}\label{eq:T-ratio}
 1+b\sum_{q\geq0}H_qz^{-q-1}
 =\frac{\cT(z)}{\cT(z+b)}.
\end{equation}

The logarithm of \eqref{eq:Cartan-charge} is triangular in the anisotropic
content moments
$C_r(\lambda)=\sum_{\square\in\lambda}c_\alpha(\square)^r$.
Equivalently, the eigenvalues satisfy
\[
 F_0=1,\quad F_1=0,\quad
 F_q=p_q^*+P_q(p_2^*,\ldots,p_{q-1}^*)\quad(q\geq2),
\]
where the $p_m^*$ are the shifted Jack power-sum coordinates and $P_q$ is a
polynomial.  Thus $H_2,H_3,\ldots$ generate the shifted Jack algebra.  Every
$\Delta_\mu^\perp$ is a polynomial in finitely many $H_q$, giving a
differential realization independent of the Nazarov--Sklyanin construction.

\subsection{Positivity in the Heisenberg basis}

Use the path operators $J_{-r}=p_r$, $J_r=\alpha D_r$.  A path for $A_q$
has heights
\[
 h_0=-1,\qquad h_{q+1}=0,\qquad h_1,\ldots,h_q\geq0.
\]
For step increments $\delta_i=h_i-h_{i-1}$ assign weight $J_{-\delta_i}$
if $\delta_i\ne0$ and weight $b h_{i-1}$ if $\delta_i=0$.  Then
\begin{equation}\label{eq:A-path}
 A_q=\sum_\gamma w_1(\gamma)\cdots w_{q+1}(\gamma).
\end{equation}

\begin{proposition}[Normally ordered positivity]\label{prop:PBW}
For every $q\geq0$, the normally ordered expansions of $A_q$ and $H_q$
have coefficients in $\NN[b]$.
\end{proposition}

\begin{proof}
In each path word, normal ordering only uses
\[
 [\alpha D_r,p_s]=\alpha r\delta_{rs}=r(1+b)\delta_{rs}.
\]
Every contraction has nonnegative coefficients in $b$, and horizontal
steps contribute $bh$ with $h\geq0$.  Hence each path has a positive
Wick expansion.  Finally, if $T=\sum t_{I,J}p_ID_J$ is normally ordered,
\[
 [D_1,T]=\sum_{I,J}m_1(I)t_{I,J}p_{I\setminus1}D_J,
\]
which preserves nonnegativity.  Apply this to $H_q=[D_1,A_q]$.
\end{proof}

\subsection{Positivity in the stable-character basis}
For a partition $\mu$ and $a\geq0$, let $\mathcal R_{\mu,a}$ be the set
of ordered sequences of horizontal-free ribbon first-return paths with
$a$ pairs and unpaired downward-step type $\mu$.  Let $\omega(R)$ be the
product of the degrees of its pairs.

For such a skeleton $R$, distinguish the internal vertices $V_1(R)$ of
the first block from all later internal vertices $V_>(R)$.  Define
$\Phi_k(R)$ by
\begin{equation}\label{eq:Phi-insertion}
 \sum_{k\geq0}\Phi_k(R)x^k
 =\prod_{v\in V_>(R)}\frac1{1-h(v)x}\frac1x
 \left(
 \frac{1+x}{\prod_{v\in V_1(R)}(1-h(v)x)}
 -\frac1{\prod_{v\in V_1(R)}(1-(h(v)-1)x)}
 \right).
\end{equation}

\begin{lemma}[Positive insertion]\label{lem:insertion}
After setting $h(v)=1+x_v$, every coefficient
$\Phi_k(R)$ belongs to $\NN[x_v:v\in V(R)]$.
\end{lemma}

\begin{proof}
Write $V_1=\{v_1,\ldots,v_s\}$ and let
$A_j=h_j(1+x_{v_1},\ldots,1+x_{v_s})$,
$B_j=h_j(x_{v_1},\ldots,x_{v_s})$, where $h_j$ is the complete homogeneous
polynomial.  The expression in parentheses in \eqref{eq:Phi-insertion},
divided by $x$, is
\[
 \sum_{j\geq0}(A_{j+1}-B_{j+1}+A_j)x^j.
\]
Every coefficient is nonnegative: $A_j$ has nonnegative coefficients,
and $A_{j+1}-B_{j+1}$ is obtained by expanding each factor $1+x_v$ and
removing the term in which no $1$ was selected.  The remaining factor
$\prod_{v\in V_>}(1-(1+x_v)x)^{-1}$ also has nonnegative coefficients.
Their product proves the claim.
\end{proof}

\begin{theorem}[Positive Cartan formula]\label{thm:positive-Cartan}
Put $d=q-|\mu|-\ellp(\mu)$.  If $d<0$, then
$[\Delta_\mu^\perp]H_q=0$.  If $d\geq0$, then
\begin{equation}\label{eq:positive-Cartan}
 [\Delta_\mu^\perp]H_q
 =\sum_{a=0}^{\lfloor d/2\rfloor}
 \alpha^{\ellp(\mu)+a}b^{d-2a}
 \sum_{R\in\mathcal R_{\mu,a}}
 \omega(R)\Phi_{d-2a}(R).
\end{equation}
In particular, $[\Delta_\mu^\perp]H_q\in\NN[b]$.
\end{theorem}

\begin{proof}
The transfer quotient \eqref{eq:T-ratio} and the Boolean definition give
\begin{equation}\label{eq:H-boolean-series}
 \sum_{q\geq0}H_qx^{q+1}
 =\frac{x+\displaystyle\sum_{n\geq2}\cB_n^{(\alpha)}x^n
       \frac{1-(1+bx)^{1-n}}b}
        {1-\displaystyle\sum_{n\geq2}\cB_n^{(\alpha)}x^n}.
\end{equation}
Expand the denominator as an ordered sequence of Boolean blocks and each
$\cB_n$ by the ribbon formula \eqref{eq:B-ribbon}.  Fix the resulting
horizontal-free skeleton $R$.  Inserting $m$ horizontal steps at a vertex
of height $h$ contributes $h^m$, hence the geometric factor
$(1-hx)^{-1}$.  The numerator singles out the first block.  Its finite
difference term is
\[
 \frac{1-(1+bx)^{1-n}}b
 =\sum_{r\geq1}(-1)^{r+1}
   \binom{n+r-2}{r}b^{r-1}x^r.
\]
The signed terms and the separate $x$ term combine, for the distinguished
block, into the bracket in \eqref{eq:Phi-insertion}.  The other blocks give
the product over $V_>(R)$.  Thus the total contribution of $R$ is
$\Phi_k(R)$, where the remaining defect is
$k=q-|\mu|-\ellp(\mu)-2a$.
Each unpaired downward step contributes $\alpha$, and each pair contributes
$\alpha$ times its degree.  The leftover $k$ horizontal units contribute
$b^k$.  Summing over $a$ and $R$ gives \eqref{eq:positive-Cartan}.
The insertion lemma, together with $\alpha=1+b$, proves positivity.
\end{proof}

\begin{example}
The first nontrivial stable-basis expansion is
\[
 H_3=3\alpha\Delta_2^\perp+2\alpha b\Delta_1^\perp.
\]
The first term comes from the horizontal-free skeleton $(+1,+1,-2)$; the
second comes from $(+1,-1)$ with one unit of remaining defect.  This
illustrates how the insertion correction produces a positive $b$ term.
\end{example}

\subsection{Polynomiality from the Cartan hierarchy}
\label{subsec:polynomiality-Cartan}

The construction of the Cartan modes by iterated commutators
yields an operator proof that the stable structure
constants are polynomials in $b$.  
For $q\geq2$, put
\[
 \psi_q=\frac{1}{\alpha}H_q.
\]
The path formula shows that $\psi_q$ has a normally ordered expansion over
$\mathbb Q[b]$.  Indeed, every nonzero path contributing to $A_q$ has a
downward step, whose weight is $J_r=\alpha D_r$ for some $r>0$.
Normal-ordering contractions use
\[
 [\alpha D_r,p_s]=\alpha r\delta_{rs},
\]
so every term retains a factor $\alpha$ before division.  After division,
all remaining weights and contractions are polynomials in $b$.

The support condition in Theorem \ref{thm:positive-Cartan} gives
no identity term, since $H_q$ annihilates the empty partition.  Therefore
\begin{equation}\label{eq:Cartan-reconstruction-triangular}
 \psi_q=q\Delta_{q-1}^{\perp}
   +\sum_{\substack{\rho\ne\varnothing,(q-1)\\
                    |\rho|+\ellp(\rho)\le q}}
       c_{q,\rho}(b)\Delta_\rho^{\perp},
 \qquad c_{q,\rho}(b)\in\mathbb Q[b].
\end{equation}
Here the coefficient of $\Delta_{q-1}^{\perp}$ is $q$: its coefficient
in $H_q$ is $q\alpha$.  Every other index in the sum has
$|\rho|\le q-2$.  Thus induction on $q$ reconstructs each
$\Delta_{(q-1)}^{\perp}$ as a normally ordered operator over $\mathbb Q[b]$.
Hall adjunction exchanges $p_r$ and $D_r$, so the same is true of
$\Delta_{(q-1)}$.

To obtain the operators with several parts, use the top-degree term of the
Jack characters in the normalization of \eqref{eq:Jack-d}:
\begin{equation}\label{eq:Jack-character-top-term}
 \operatorname{gr}_{|\mu|}\Ch_\mu^{(\alpha)}
   =\alpha^{|\mu|-\ellp(\mu)}p_\mu,
 \qquad
 p_\mu(x)=\prod_i\sum_jx_j^{\mu_i}.
\end{equation}
This is the standard top-degree characterization of Jack characters
\cite{BenDaliDolega}.  If $\ellp(\mu)\ge2$, all parts of $\mu$ are
strictly smaller than $|\mu|$.  Assuming inductively that the corresponding
single-part operators have been reconstructed, the composition
\[
 T_\mu=\Delta_{\mu_1}\cdots\Delta_{\mu_{\ellp(\mu)}}
\]
has a normally ordered expansion over $\mathbb Q[b]$.  Its eigenvalue has
the same top-degree term as $\Ch_\mu^{(\alpha)}$.  Since the power sums
$p_\rho$ of a fixed degree are linearly independent, expansion in the
Jack-character basis therefore has the triangular form
\[
 T_\mu=\Delta_\mu+
       \sum_{|\rho|<|\mu|}e_{\mu,\rho}(b)\Delta_\rho.
\]
The balanced-action lemma shows that the coefficients
$e_{\mu,\rho}(b)$ lie in $\mathbb Q[b]$: apply $T_\mu$ to the unit
$\widehat p_{\varnothing}=1/(1-p_1)$ and use its normally ordered form.
Subtracting the already constructed lower-degree operators now gives
$\Delta_\mu$ over $\mathbb Q[b]$.  This completes the induction for every
partition $\mu$.

Finally, applying this normally ordered form of $\Delta_\mu$ to
$\widehat p_\nu$ and using Lemma \ref{lem:balanced-action} gives
\[
 d_{\mu\nu}^{\lambda}(1+b)\in\mathbb Q[b]
 \qquad\text{for all }\lambda,\mu,\nu.
\]
This recovers the polynomiality of the stable structure constants by a
different route.  The argument uses the positive Cartan expansion and the
top-degree Jack-character property; its subtractions do not preserve
positivity, so it gives no new positivity statement.

For example, the first steps are
\[
 \psi_2=2\Delta_1^\perp,\qquad
 \psi_3=3\Delta_2^\perp+2b\Delta_1^\perp.
\]
Since $\Ch_1(\lambda)=|\lambda|$ and
$\Ch_{11}(\lambda)=|\lambda|(|\lambda|-1)$, one has
$(\Delta_1^\perp)^2=\Delta_{11}^\perp+\Delta_1^\perp$.  Using the
displayed formula for $H_4$ in the exact computations, one obtains
\[
 \Delta_3^\perp
 =\frac14\left(
 \psi_4-\frac32\alpha \psi_2^2-3b\psi_3
       +2(1+b+b^2)\psi_2\right),
 \qquad \alpha=1+b.
\]
This illustrates how the induction recovers a stable multiplier using
only Cartan modes and products of lower modes.

\subsection{Integrality of the normally ordered coefficients}
\label{subsec:integrality}

The integral normalization of the Boolean generators allows one to
replace $\QQ[b]$ by $\mathbb Z[b]$ in the polynomiality statement of the
previous subsection.

\begin{theorem}\label{thm:normal-integrality}
For every partition $\mu$, the normally ordered coefficients of
$\Delta_\mu^\perp$ and $\Delta_\mu$ belong to $\mathbb Z[b]$.
Moreover,
\[
 \alpha^{\ellp(\mu)-\ellp(J)}
 [p_ID_J]\Delta_\mu^\perp\in\mathbb Z[b]
 \qquad\text{for every }I,J.
\]
In particular, the stable structure constants in \eqref{eq:Jack-d}
belong to $\mathbb Z[b]$.
\end{theorem}

\begin{proof}
We work first over $\QQ(\alpha)$. In the normalization of this paper,
for $N=|\lambda|\geq |\mu|$ the stable character is
\begin{equation}\label{eq:integral-character-normalization}
 \Ch_\mu^{(\alpha)}(\lambda)
 =z_\mu\binom{N-|\mu|+m_1(\mu)}{m_1(\mu)}
 [p_{\mu\cup1^{N-|\mu|}}]J_\lambda^{(\alpha)}.
\end{equation}
For $N<|\mu|$ it is zero. To see the normalization, the degree-$N$
component of $\whp_\mu$ is
$N!p_\mu p_1^{N-|\mu|}/(N-|\mu|)!$. The Hall pairing of this component
with $J_\lambda^{(\alpha)}$, divided by $N!$ as in
\eqref{eq:Jack-product}, gives \eqref{eq:integral-character-normalization}
after writing out $z_{\mu\cup1^{N-|\mu|}}$. These characters are the
eigenvalues of $\Delta_\mu^\perp$ on $J_\lambda^{(\alpha)}$.
Indeed, the vectors $J_\lambda^{(\alpha)}[\alpha X]$ and
$J_\lambda^{(\alpha)}[X]$ are dual up to their Jack norms under the Hall
pairing, by its definition on power sums and Jack orthogonality.

We use the $\mathbb Z[\alpha]$ version of the integral generation
theorem of Ben Dali and Do\l{}\k{e}ga
\cite[Theorem~4.2.12]{BenDaliThesis}. In the normalization
\eqref{eq:integral-character-normalization}, it identifies the
$\mathbb Z[\alpha]$-algebra of stable characters with the algebra
generated by the Boolean cumulants divided by $\alpha$. Thus each
$\Delta_\mu^\perp$ is a finite polynomial over $\mathbb Z[\alpha]$ in
$\cB_m^{(\alpha)}/\alpha$, $m\geq2$. Here the eigenvalues of the
operator generators agree with the Boolean cumulants by the
Nazarov--Sklyanin identification in Section~4. Equality on the Jack
basis verifies the operator identity on each homogeneous component.

By Proposition~\ref{prop:first-return},
\[
 \frac{\cB_m^{(\alpha)}}{\alpha}
 = L_{0+}L_{++}^{m-2}(D_1,D_2,\ldots)^{\mathsf t}.
\]
The division removes the common factor in
$L_{+0}=(\alpha D_1,\alpha D_2,\ldots)^{\mathsf t}$ before any normal
ordering is performed. Every remaining matrix entry in
\eqref{eq:Lax} has coefficients in $\mathbb Z[b]$. In moving derivatives
past multiplications, the only additional scalar factors are the integers
from $D_i p_j=p_jD_i+i\delta_{ij}$. Each generator and each finite product
of generators therefore has normally ordered coefficients in
$\mathbb Z[b]$.

For the local finiteness in this argument, give $f w^i$ the total
weight $\deg(f)+i$. Each entry $L_{ij}$ maps the $j$th auxiliary
component to the $i$th and preserves this total weight. At weight $N$
there are only the heights $0,\ldots,N$ to consider. Products of the
resulting operators preserve weight, and their normally ordered terms
with $|J|>N$ act as zero on weights at most $N$. The coefficients
with $|I|=|J|\leq N$ are determined by this truncated action. For
example, one can apply the operator in $X$ to the Hall Cauchy kernel,
divide by that kernel, and read the coefficients of $t^k p_I[X]p_J[Y]$
for $k\leq N$. This also shows that computing with finite truncations
gives consistent coefficients as $N$ increases. Since
$\mathbb Z[\alpha]=\mathbb Z[b]$, this proves integrality for
$\Delta_\mu^\perp$. Hall adjunction sends $p_ID_J$ to $p_JD_I$ and
proves integrality for $\Delta_\mu$.

For the normalized coefficients, let
\[
 \tau_\alpha(t;A,X,Y)=\sum_\lambda t^{|\lambda|}
 \frac{J_\lambda^{(\alpha)}[A]J_\lambda^{(\alpha)}[X]
       J_\lambda^{(\alpha)}[Y]}{j_\lambda^{(\alpha)}}.
\]
Set $p_r[\varepsilon]=\delta_{r1}$ and
\[
 K_\alpha(t;X,Y)=\tau_\alpha(t;\varepsilon,X,Y)
 =\exp\!\left(\sum_{r\geq1}\frac{t^r p_r[X]p_r[Y]}{\alpha r}\right).
\]
Let $U_\mu(t;X,Y)$ denote the sum of the normalized coefficients in the
theorem, with $p_ID_J$ replaced by $t^{|J|}p_I[X]p_J[Y]$. Applying
$\Delta_\mu^\perp$ to $K_\alpha$ in $X$ and dividing by $K_\alpha$
replaces $D_j$ by $t^j p_j[Y]/\alpha$. For the alphabet
$\varepsilon+\alpha Z$, defined by
$p_r[\varepsilon+\alpha Z]=\delta_{r1}+\alpha p_r[Z]$, we obtain
\begin{equation}\label{eq:normalized-symbol-kernel}
 \sum_\mu\frac{p_\mu[Z]}{z_\mu}U_\mu(t;X,Y)
 =\frac{\tau_\alpha(t;\varepsilon+\alpha Z,X,Y)}
 {K_\alpha(t;X,Y)}.
\end{equation}
Indeed, the coefficient of $p_\mu[Z]$ in
$J_\lambda^{(\alpha)}[\varepsilon+\alpha Z]$ is
$\alpha^{\ellp(\mu)}\Ch_\mu^{(\alpha)}(\lambda)/z_\mu$.
The binomial coefficient in \eqref{eq:integral-character-normalization}
counts the choices of the $m_1(\mu)$ occurrences of $p_1[Z]$.

The coefficients of $\alpha\log\tau_\alpha$ are polynomials in
$\alpha$ over $\QQ$ by \cite[Theorem~1.2]{DolegaFerayBCumulants}.
Form the difference of the two specializations with
$A=\varepsilon+\alpha Z$ and $A=\varepsilon$. The terms independent
of $Z$ cancel. Every term that remains contains at least one
$\alpha p_r[Z]$ and is therefore divisible by $\alpha$ over
$\QQ[\alpha]$. After this division we obtain the logarithm of the
right-hand side of \eqref{eq:normalized-symbol-kernel}. Its constant
term in $t$ is zero, so exponentiation and extraction of
$p_\mu[Z]$ show that the coefficients of $U_\mu$ are in
$\QQ[\alpha]$. By the first part of the proof they are also in
$\mathbb Z[\alpha,\alpha^{-1}]$. The intersection of these rings
inside $\QQ(\alpha)$ is $\mathbb Z[\alpha]$, which proves the
stronger normalization assertion.

For the stable structure constants, apply
Lemma~\ref{lem:balanced-action} to each integral normal monomial.
For fixed $\nu$ and a given output index $\lambda$, only finitely many
monomials contribute. Indeed, the input derivative indices at least
two must be contained in $\nu$, and \eqref{eq:balanced-action} bounds
$j_1$ by $|\lambda|$. Each contribution is in $\mathbb Z[b]$.
The empty character gives the identity operator, for which the
assertions hold as well.
\end{proof}

\section*{Acknowledgments}
The author used ChatGPT (OpenAI) as an interactive aid in organizing
the manuscript, editing the English text, bibliographical searches, and for numerical and
computer-algebra calculations.  All mathematical statements and
conclusions remain the responsibility of the author.

\end{document}